\documentstyle[a4,twoside]{article}
\makeatletter
\renewcommand{\@oddhead}{
Nonlinear Approximation of an Operator Equation 
\hfill \thepage}
\renewcommand{\@evenhead}{\thepage \hfill 
Nonlinear Approximation of an Operator Equation 
}
\renewcommand{\@oddfoot}{}
\renewcommand{\@evenfoot}{}
\makeatother
{\vspace*{2ex}\par {#1}\quad\it }{\medskip\par }

\newenvironment{Thm}[2]{\par\addvspace{\bigskipamount}{\bf #1#2}\it }%
{\par\addvspace{\bigskipamount} }

\newenvironment{Definition}[1]{\begin{Thm}{Definition}{#1}}{\end{Thm} }

\newenvironment{Theorem}[1]{\begin{Thm}{Theorem}{#1}}{\end{Thm} }

\newenvironment{Proposition}[1]{\begin{Thm}{Proposition}{#1}}{\end{Thm} }

\newenvironment{Corollary}[1]{\begin{Thm}{Corollary}{#1}}{\end{Thm} }
\newenvironment{Observation}[1]{\begin{Thm}{Observation}{#1}}{\end{Thm} }
\newenvironment{Remark}[1]{\begin{Thm}{Remark}{#1}}{\end{Thm} }

\newenvironment{Example}[1]{\begin{Thm}{Example}{#1}}{\end{Thm} }

\newenvironment{Proof}{\par\addvspace{\bigskipamount} {\bf Proof}}%
{\par\hspace*{\fill}$\Box$ \par\addvspace{\bigskipamount} }

\author{S.A.~Choro\v{s}avin}
\title{A Nonlinear Approximation \\
of Operator Equation
$V^{*}QV=Q$  : \\
Nonspectral Decomposition of Nonnormal Operator
and Theory of Stability}

\date{}
\newcommand{\R}{Ran \,Y_0^{1/2}}
\newcommand{\E}{VY_0V^*=\frac{Y_0}{I+Y_0}}
\newcommand{\bra}[1]{\langle{#1}|}

\begin{document}
\maketitle

\begin{abstract}
    $V$ denotes arbitrary bounded bijection on Hilbert space $H$.
We try to describe the sets of $V$-stable vectors,
i.e.
$\{x\in H|
\mbox{ the sequence }\|V^N x\| (N=1,2,\ldots)\mbox{ is bounded}\}$
and some other analogous sets.
We do it in terms of one-parameter operator equation
$ Q=V^*(Q+tI)(I+tQ)^{-1}V $
($t$ is real valued parameter $0\leq t \leq 1$, 
$Q$ is operator to be found, $0\leq Q $)
\end{abstract}
%\vspace{8mm}
 
\section{Introduction}
 
   Throughout this paper $H$ will denote a Hilbert space with scalar product 
$<,>$, $V$ denotes a linear bounded bijection $H$ onto $H$,
$$
    r(T):= \hbox{spectral radious of }T \;.
$$
We will discuss the structure of the next four sets:

%\begin{Definition}{ D1-1}
%\end{Definition}
%
$$
Stab_+(V)\;:=\;\{x\in H|\forall a>1 \; \exists C\geq 0 \; \forall N\geq 0 
\quad\|V^Nx\|\leq Ca^N\}.
$$
$$
 Stab (V) :=\{x\in H|\quad \|V^nx\| \leq M \mbox{ for some real } M
 \mbox{ and every }n=0,1,2, \dots \}
$$ 
$$
 Stab_0 (V) :=\{x\in H|\quad \|V^nx\| \to 0\quad (n\to \infty )\}
$$ 
$$ 
 l_2 (V)
:=\{x\in H|\quad \|Vx\|^2+\|V^2x\|^2+\cdots+\|V^nx\|^2+\cdots<\infty \} 
$$

    Recall $V$ is 
{\bf similar to an unitary operator} iff there exists a bounded
uniformly positive operator $Q$ such that
$$
         V^*QV=Q
$$
With this equation we shall consider an 'approximation' equation
( parametrized by real $t$)
$$
  Q=V^*\frac{Q+t}{I+tQ} V,\quad Q\geq 0,\quad 0<t\leq 1 \qquad \eqno(^ *)
$$
(hereinafter $t$ denotes always a real number such that $0<t\leq 1$
and if no confusion can occur we shall often write $t$ instead of
$tI$, $I$ is identity operator).

     The interest in this equation can be motivated by the next

\begin{Example}{ 1.} Let $V$ be normal. Routine, though tedious
calculation shows that 
$$
Q_t:=\left[ -\frac{I-V^*V}{2}+
\sqrt{\left( \frac{I-V^*V}{2}\right) ^2+t^2V^*V}\,\right]\frac{1}{t}
$$
is uniformly positive solution of $(^*)$, there holds  
$$ 
Q_t^{-1}=\left[ -\frac{I-(V^*V)^{-1}}{2}+
\sqrt{\left( \frac{I-(V^*V)^{-1}}{2}\right) ^2+t^2(V^*V)^{-1}}\,\right]
\frac{1}{t}  
$$
and there exist
$$
 X_0\;:=\;strong-lim\;tQ_t =(V^*V-I)E(1,\infty ) \qquad (t\to +0)
$$
$$
 Y_0\;:=\;strong-lim\;tQ_t^{-1}=( ( V^*V)^{-1}-1)E(0,1) \qquad (t\to +0)
$$
Besides that it is fast evident that $Q_tE([1])=E([1]),\;
Q_tE(1,\infty ) $ is monotone increasing (with $t\to +0$),
$Q_tE(0,1)$ is monotone decreasing and there exists an
$$
R_0\;:=strong-lim(I+Q_t)^{-1} = E(0,1)+E([1])/2 \qquad (t\to +0).
$$
here $E(\Delta )$ denotes the spectral function of the selfadjoint
operator $V^*V$.
Note:
$$
I-R_0 = E([1])/2+E(1,\infty )
$$
$$
Ker(I-R_0)=E(0,1)H
$$
$$
RanX_0 \subset E(1,\infty )H =\overline{RanX_0} 
$$
$$
RanY_0 \subset E(0,1)H=\overline{RanY_0}
$$

                           $\Box $
\end{Example}

    Also, in the above considered case of the normal $V$ it is 
estableshed that the operators $X_0,Y_0,R_0$ define (in essential) 
the spectral subspaces of $V$ (with $V$ together one can consider
$aV-b, b/a \not\in   spectrum\, V$) . In this article we shall show
that the similar situation holds for the arbitrary bounded bijection $V$.
\medskip

We follow standards of [RS] when we apply  mathematical concepts
and sometimes we apply P.A.M. Dirac's `bra-ket' syntax.

\newpage
  We will often  cite some assertions and propositions of [Ch1-4].
For the most convenient and accesible 
way to do it, we collect them together and  resume them here as

\begin{Theorem}{ 1.}

{\rm(i)} the solution of $ (^*)$ exists and it is unique; denote it by $Q_t$

{\rm(ii)} $Q_t$ is bounded selfadjoint uniformly positive and there are 
satisfied
     inequalities:
$$
   tV^*V\leq Q_t \leq V^*V/t
$$
$$
   tQ_t\leq sQ_s \qquad (0<t\leq s\leq 1)
$$
denote $X_t:=tQ_t $

{\rm(iii)} $Q_t ^{-1}$ is (unique) solution of the analogous equation:
$$
  Q_t ^{-1} =V^{-1} \frac{Q_t^{-1}+t}{I+tQ_t^{-1}} (V^{-1})^* 
$$
so there are satisfied inequalities:
$$
  tV^{-1}(V^{-1})^* \leq Q_t^{-1} \leq V^{-1}(V^{-1})^*/t
$$
$$
   tQ_t^{-1}\leq sQ^{-1}_s  \qquad (0<t\leq s\leq 1)
$$
denote $   Y_t:= tQ_t^{-1}   $.

{\rm(iv)}  Let
$$
  X_0:=strong-lim\;X_t,\quad Y_0:=strong-lim\;Y_t  \quad (t\to +0)
$$
The operators $X_0,\; Y_0$ are bounded positive and they are maximal
solutions of the equations
$$
X=V^*\frac{X}{I+X}V,\; X\geq 0,\;\mbox{resp.}\;VYV^*=\frac{Y}{I+Y},\;Y\geq 0
$$
('maximal' denotes 'maximal with respect to usual partial order for
bounded operators on Hilbert space')

{\rm(v)} There hold the formulae
$$
 Y_0=strong-lim\;(V^*V+V^{*2}V^2+\cdots+V^{*n}V^n)^{-1},\;(n\to \infty )
$$
$$
X_0=strong-lim\;((V^*V)^{-1}+(V^{*2}V^2)^{-1}+\cdots+(V^{*n}V^n)^{-1})^{-1},
                          \;(n\to \infty )
$$

{\rm(vi)}  Denote  $ R_t:=(I+Q_t)^{-1}  $. Then $0\leq R_t\leq I$,
     $Q_t=R_t^{-1}-I$ and the equation $(^*)$ is equivalent to the
     equation
$$
[R_t+t(I-R_t)](V^{-1})^*(I-R_t)=[(I-R_t)+tR_t]VR_t
$$

{\rm(vii)} Let $R_0$ be a weak operator limit point of the net 
$\{R_t,\,t\to +0\}$
(it is clear that $R_0$ exists and $0 \le R_0 \le I$). Then
$$
V\,Ker\,(I-R_0)=Ker\,(I-R_0) %%%%%%%%%%%%\qquad \Box
$$
$$
V^{*-1}\,Ker\,R_0=Ker\,R_0   %%%%%%%%%%%%%%\qquad \Box
$$
In particular,
$$
V\overline{Ran\,R_0}=\overline{Ran\,R_0} \qquad \Box
$$

\end{Theorem}

\newpage

\section{ Equation $\displaystyle   Q=V^*\frac{Q+t}{I+tQ}V$ . 
          General Properties.}

Hereinafter {\bf F} denotes an arbitrary ultra filter, which
majorizes usual convergence to $+0$.
We will write 
$$
  t \to +0 \quad \mbox{ instead of }  t\stackrel{\bf F}{\longrightarrow} +0
$$
if no confusion can occur.
\begin{Definition}{ 1.}
$$ Fin\,Q:=\{x\in {\bf H} | <x,Q_tx>\leq M_x \mbox{ for some real }M_x
    \mbox{ and for almost every }  t
       \mbox{ resp. }{\bf F } \}
$$

$$ Ker\,Q_0:=\{x\in {\bf H} | <x,Q_tx>\to 0 \mbox{ for } t\to +0 \}
                                %%%%%%%%%%%%%%\mbox{ with respect to }{\bf F } 
$$
\end{Definition}

\begin{Theorem}{ 1.}
$$
  Ran\,Y_0^{1/2}=l_2(V)
      \subset Stab_0(V) 
      \subset Ker\,Q_0 
      \subset Ker\,(I-R_0) 
      \subset \overline{ Ran\,R_0 } 
      \subset Ker\,X_0    \eqno(*)
$$
$$
  Stab(V) \subset Fin\,Q \subset \overline{Ran\,R_0}\subset Ker\,X_0  \eqno(**)
$$
Every set of these series is $V$-surinvariant.
(recall, some {\bf $L$ is said to be $T$-surinvariant, iff $TL=L$}).

In addition
$$
     Fin\,Q \subset { Ran\,R_0^{1/2} } 
$$
\end{Theorem}

\begin{Observation}{ 0.} 
It is well-known and evident that
$$
\begin{array}{cccccc}
 A^*A\leq B^*B 
&\iff&
\|Ax\|\,\le\,\|Bx\|\quad(x\in H)
&\Longrightarrow&
RanA^*\;\subset \;RanB^*\\[5pt]
 & & & \makebox[0ex][c]{($A$ and $B$ are bounded)}. 
\end{array}
$$
\end{Observation}

\begin{Corollary}{ 0.} 
Let 
$x\in H$ , $Y$ be  selfadjoint and let  $Y\geq 0$. Then
$$
 x \in Ran\,Y^{1/2}  \iff x\bra{x}\leq cY \qquad  \mbox{ for some real $c$. } 
$$
\end{Corollary}
\begin{Proof}{.}

Proof of
$
x\bra{x}\leq cY 
 \Longrightarrow x \in Ran\,Y^{1/2}    
$:

Take into account Observation 0. Then obtain
$$
x\bra{x} \leq cY  
\Longrightarrow Ran\,x\bra{x} \subset Ran\,Y^{1/2}  
\Longrightarrow x \in Ran\,Y^{1/2}   \,.
$$

Now proof of 
  $x \in Ran\,Y^{1/2} \Longrightarrow  x\bra{x} \leq cY$ :

By the definition of $Ran$ 
$$
x \in Ran\,Y^{1/2}  
\iff 
x=Y^{1/2}y    \quad \mbox{ for an $y$.}
$$
Hence 
$$
x \in Ran\,Y^{1/2}  
\Longrightarrow 
x\bra{x} 
\equiv 
Y^{1/2}y\bra{Y^{1/2}y }
=Y^{1/2}y\bra{y}Y^{1/2}
\leq \|y\|^2 Y        \,.
$$
Now denote $\|y\|^2$  by $c$.

\end{Proof}

\newpage

\begin{Observation}{ 1.} 
Let $\{A_t\}_t$ be a net of selfadjoint 
positive bounded operators. Suppose $A_t \le aI$ for some positive number
$a$ (and every $t$), $A_0$ be a weak operator limit point of this
net (clear: $A_0$ exists and $0 \le A_0 \le aI$ ). 

Then
$$
<x,A_t x>\to 0\Leftrightarrow \|A_tx\|\to 0\Leftrightarrow x\in Ker\;A_0 
$$
\end{Observation}
\begin{Proof}. 
Clear (see e.g. [Ch1]).
\end{Proof}

%\newpage

\begin{Observation}{ 2.} Let $Q$ be bounded selfadjoint positive,$0<t\le 1$ .
Then
$$
\frac{Q}{I+t\|Q\|} \le \frac{Q+t}{I+tQ}\le Q+t
$$
\end{Observation}
\begin{Corollary}{.}
$$\fbox{$\displaystyle
VFin\,Q=Fin\,Q\,,\quad  VKerQ_0=KerQ_0 \,.
$}$$
\end{Corollary}
\begin{Observation}{ 2'.} Let $Q$ be bounded selfadjoint positive,$0<t\le 1$ .
Then
$$
\frac{1-t}{1+t} Q +\frac{2t}{1+t} -\frac{Q+t}{I+tQ} =
                   \frac{t(1-t)}{1+t}  \frac{(Q-I)^2}{I+tQ}
$$
In particular
$$
\frac{Q+t}{I+tQ} \le \frac{1-t}{1+t} Q +\frac{2t}{1+t} 
$$
Denote $ z:=(1-t)/(1+t)$. Then $ 1-z=2t/(1+t) $ and with these denotations
$$ 
\frac{Q+t}{I+tQ} \le zQ+(1-z)
$$
Note
$$
t\to +0 \iff z\to 1-0 \,.
$$
\end{Observation}

\begin{Proof}. Clear.
\end{Proof}

%\newpage

\begin{Observation}{ 3} . For $Q_t$ the just mentioned inequality gives
$$
   Q_t\; \le \;V^*(zQ_t+(1-z))V
$$
and with iterating this inequality one can obtain
$$
Q_t\le (1-z)[V^*V+zV^{*2}V^2+\cdots+z^{n-1}V^{*n}V^n]+z^nV^{*n}Q_tV^n
$$
In particular, given numbers $z, M_0$ and an $x\in H $ such that
$$
0 <z<1,\; M_0\ge 0,\; \|V^nx\|^2\le M_0\;(n=1,2,3,...)
$$
then
$$
<x,Q_tx>\le (1-z)(\|Vx\|^2+z\|V^2x\|^2+z^2\|V^3x\|^2+\cdots)\le M_0 \,.
$$

In particular

$$\fbox{$\displaystyle 
Stab(V) \subset Fin\,Q
$}$$

Now suppose $\|V^nx\| \to 0$ ($n\to \infty$) for
some $x\in H$ .
Let $t \to +0$. Then $z \to 1-0$ and hence $<\!x,Q_tx\!>\;\to 0$
since $\|V^nx\|\to 0$ and  $<x,Q_tx>\;\ge 0$. 

  In particular
$$\fbox{$\displaystyle 
Stab_0(V)\,
\subset \,
Ker Q_0
$}$$
\hspace*{\fill} $\Box$
\end{Observation}

\newpage

\begin{Observation}{ 4.}
Let
$x \in H$, $\|x\|=1$. Then

a)
$$
x\in Fin\,Q \Rightarrow \exists M\geq 0 \forall \alpha \geq 0
: \qquad
x\bra{x} \leq (M+\alpha)(Q_t+\alpha)^{-1}
$$

b) Given some $M, \alpha \geq 0$,
such that
$$
x\bra{x} \leq (M+\alpha)(Q_t+\alpha)^{-1}
$$
then
$$
x\in Fin\,Q
$$
\end{Observation}

%%% Insbesondere,
\begin{Proof}
\begin{eqnarray*}
<x, Q_tx>\leq M &&\\
&\iff&  <x,(Q_t+\alpha)x> \leq M+\alpha \\ 
&\iff&  \|(Q_t+\alpha)^{1/2}x\|^2 \leq M+\alpha  \\
&\iff&  (Q_t+\alpha)^{1/2}x\bra{(Q_t+\alpha)^{1/2}x} \leq M+\alpha \\ 
&\iff&  (Q_t+\alpha)^{1/2}x\bra{x}(Q_t+\alpha)^{1/2} \leq M+\alpha \\ 
&\iff& x\bra{x} \leq (M+\alpha)(Q_t+\alpha)^{-1} 
\\ 
\end{eqnarray*}
\end{Proof}

\begin{Corollary}{.}
\fbox{$\displaystyle
Fin\,Q \subset Ran\, R_0^{1/2} \subset \overline{Ran\,R_0}
$}
\end{Corollary}

\begin{Proof}
First apply Observation 4 for 
$\alpha:=1$
and let 
$t\to +0$.
Then obtain
$$
 x\bra{x} \leq (M+1)R_0;
$$ 
Now apply Observation 0 or Corollary 0.
\end{Proof}

%\newpage

\begin{Observation}{ 5.}
Let 
$x\in H$ . Then  %%%%%Dann gilt
$$
x \in \R  \iff x\bra{x}\leq cY_0 \iff <x, Q_t x> \leq ct
$$
\end{Observation}

\begin{Proof}
$$
x \in \R  
\iff 
x=Y_0^{1/2}y 
\Longrightarrow 
x\bra{x} 
\equiv 
Y_0^{1/2}y\bra{Y_0^{1/2}y }
\leq \|y\|^2 Y_0 
$$
Besides,    %%%%%%%%Au\ss erdem
$$
x\bra{x} \leq cY_0  
\Longrightarrow Ran\,x\bra{x} \subset \R  
\Longrightarrow x \in \R 
$$
Now recall that    %%%%%%%Nun erinnern wir daran, da\ss
$$
Y_0 \leq Y_t \mbox{ und  } Y_t \stackrel{s}{\to} Y_0 \,.
$$
Hence                    %%%%%%%Somit ergibt sich, da\ss
\begin{eqnarray*}
x \in \R  
&\iff& x\bra{x}\leq cY_0 \\
&\iff& x\bra{x}\leq cY_t \\
&\iff& x\bra{x}\leq ctQ_t^{-1} \\
&\iff& Q_t^{1/2}x\bra{Q_t^{1/2}x}\leq ct \\
&\iff& \|Q_t^{1/2}x\|^{1/2}\leq ct \\
&\iff& <x, Q_t x> \leq ct \\
\end{eqnarray*}

%%%%Insbesondere,

%%%%%$$
%%%%%%Ran Y_0^{1/2} \subset Ker Q_0
%%%%%$$

\end{Proof}

%\newpage

\begin{Proof}{ \bf of } 
\fbox{$\displaystyle
Ker Q_0 \subset \,Ker (I-R_0) .
$}

Note
$$
I-R_t=Q_t(I+Q_t)^{-1}\leq Q_t
$$
Now suppose
$x\in Ker Q_0$ 
i.e. $<x, Q_tx> \to 0$ 
and apply the Observation 1 to $A_t=I-R_t$.

\end{Proof}
\begin{Proof}{ \bf of} 
\fbox{$\displaystyle
Ker (I-R_0) \subset \overline{Ran\,R_0}.
$}

Note $Ker\;$ is closed and use definitions of $Ker\;,  Ran\,$.
\end{Proof}

\begin{Proof}{ \bf of} 
\fbox{$\displaystyle
\overline{Ran\,R_0}\subset Ker\,X_0 \,.
$}

Note
$$
R_tX_t = t(I+Q_t)^{-1} = X_tR_t
$$
Then 
$$
R_tX_t = X_tR_t \stackrel{\|\cdot\|}{\longrightarrow} 0.
$$
Recall $R_t, X_t$ are selfadjoit, positive, bounded
and
$X_t \stackrel{s}{\longrightarrow} X_0$, 
$R_t \stackrel{w}{\longrightarrow} R_0$.
Hence
$R_0 X_0= 0 = X_0 R_0$.

The rest is evident.
\end{Proof}

\newpage
 
\begin{Proof}{ of} 
\fbox{$\displaystyle
Ran\,Y_0^{1/2}\; \subset \; l_2(V)\, .
$}

We have
$$
             \E          \eqno(^*)
$$
Hence $\|Y_0^{1/2}V^*x\|=\|(I+Y_0)^{1/2}Y_0^{1/2}x\|\, \le \, \|Y_0^{1/2}x\|$
and the relation
$$
W:\, Y_0^{1/2}x\, \to \, Y_0^{1/2}V^*x
$$
defines a contraction $Ran\,Y_0^{1/2}\, \to \, \R$. 
This contraction has extension
to a contraction $\overline{\R}\, \to \, \overline{\R}$. It will be denote
by $W $ too. Clear, there hold
$$
\begin{array}{ll}
(i)\; WY_0^{1/2}=Y_0^{1/2}V^* & (iii)\; W^*W=(I+Y_0)^{-1}|\overline{\R}
\\[0.5cm]
(ii)\; Y_0^{1/2}W^*P=VY_0^{1/2} & (iv)\; Y_0P=((W^*W)^{-1}-I)P
\end{array}
$$
and $ WP=PWP,\quad W^*P=PW^*P  $ .Here is
$$
  P:= orthoprojection \;  onto \;  \overline{Ran\,Y_0^{1/2}}
$$
Hence
$$
 Y_0^{1/2}V^{*n}V^nY_0^{1/2}\;=\;W^nY_0W^{*n}P\;=\;PW^nY_0W^{*n}P
$$
$$
     =\;PW^n((W^*W)^{-1} -I)W^{*n} P
$$
$$
       =\;P(W^{n-1}W^{*n-1}-W^nW^{*n})P.
$$
and
$$
    \sum_{n=2}^N \|V^nY_0^{1/2}x\|^2  
      = \|WPx\|^2  - \|W^{N}Px\|^2
$$
%But $  W^*   $  is a contraction. 
Hence
$$
    \sum_0^\infty \|V^nY_0^{1/2}x\|^2  < \infty   \qquad \Box
$$
\end{Proof}

\begin{Remark}{ to this proof.}

It follows from $(^*)$ and $  Y_0\ge 0  $ that
$$
  VRan\;Y_0\;=\;Ran\;Y_0,\quad V^*Ker\;Y_0\;=\;Ker\;Y_0,\quad 
 V\overline{Ran\;Y_0}  \;=\; \overline{Ran\;Y_0} 
$$
Besides that,
$$
\begin{array}{rcl}
(iii) & => & \|Wx\|\, \ge\, \|x\|/(1+\|Y_0\|) \quad (x\in \,D_W)\\[2mm]
(ii)  & => & V\,\R\; \subset \; \R
\end{array}
$$
$$
 \mbox{($i$, or definition of $W$)}\;=>\; W\R\,=\,Ran\,Y_0^{1/2} \qquad \Box
$$

\end{Remark}

\newpage
 
\begin{Proof}{ of } 
\fbox{$\displaystyle
l_2(V)\; \subset \; Ran\,Y_0^{1/2}      \,.
$}
 
Let $x \in l_2(V)$ and set
$$
\begin{array}{rcl}
c & := & \|Vx\|^2+\|V^2x\|^2+\|V^3x\|^2+\cdots\\[5pt]
P & := & orthoprojection\; mapping\; onto\; span\{x\}
\end{array}
$$
Then
$$
<\!x, Y_t^{-1}x\!>\,=\,<\!x, \frac{1}{t}Q_tx\!>\; \le \;
\frac{1-z}{t}\,c\; =\; \frac{2c}{1+t}\; \le \; 2c
$$
Hence
$$
\|x\|^2P\, \le \, 2cY_t,\quad \|x\|^2P\, \le \, 2cY_0,\quad
Ran\,P \; \subset \; Ran\,Y_0^{1/2}\quad 
$$
\end{Proof}

\begin{Remark}{ to this proof. }
$$
\begin{array}{rl}
(i) & Y_t\,\gg \,0,\; Y_0\,\ge \, 0, \; P=P^{1/2}\, \ge \,0 \\[5pt]
(ii)& \|Ax\|\,\le\,\|Bx\|\quad(x\in H)\;=>\;RanA^*\;\subset \;RanB^*\\[5pt]
    & \mbox{($A$ and $B$ are bounded)}. \qquad  \Box
\end{array}
$$
\end{Remark}

\begin{Corollary}{. }  
\fbox{$\displaystyle
l_2(V) = Ran\,Y_0^{1/2}  \,.
$}
\end{Corollary}
 
\newpage

\section{ Nonspectral Decomposition}
\begin{Observation}{ 1}.

  It follows from the definitions of $  X_t,\;Y_t $ that
$  X_tY_t\;=\;Y_tX_t\;=\;t^2 $ . Hence $ X_0Y_0\;=\;Y_0X_0\;=\;0 $ 
 and $ \overline{Ran\;X_0} \subset Ker\;Y_0$   ,
     $ \overline{Ran\;Y_0} \subset Ker\;X_0$  .
   But $X_0,\;Y_0$ are selfadjoint. Thus we obtain an orthogonal
decomposition
$$
H\;=\;\overline{Ran\;Y_0}\;+(Ker\;X_0\cap Ker\;Y_0)\;+
                              \overline{Ran\;X_0}
$$
such that

1) first  component is $V$-surinvariant;

2) third component is $V^{*-1}$-surinvariant.

(recall, some {\bf $L$ is said to be $T$-surinvariant, iff $TL=L$}).
      Moreover, denote $j_t:=(I+tQ_t^{-1})^{1/2}/(I+tQ_t)^{1/2}$ ,then 
  $j_t$ is uniformly positive, bounded, there hold
$$
 (1+\|V\|^2)^{-1/2}\;\le\;j_t\;\le\;(1+\|V^{-1}\|^2)^{1/2}\;,
$$
and $  (j_tV)^*Q_t(j_tV)\,=\,Q_t$. In particular $j_tV$ is similar to an
unitary operator.

  For a moment suppose
$dim\,H\,<\infty $.
It is clear that now 
the restriction of $V$ onto $\overline{Ran\,Y_0}$  and 
the restriction of $V^{*-1}$ onto $\overline{Ran\,X_0}$  
are similar to uniform contractions (see theorems 2.1 with remarks to
the proof  of $l_2(V)=\R$).
In addition, if
$$RanY_0=\{0\}=RanX_0\,,$$
then the restriction of $V$ onto $Ker\,X_0\cap Ker\,Y_0$ has the unit
spectrum.  $\Box$  
\end{Observation}

This motivates the

\begin{Definition}{ 1.}
 We shall say, a linear bounded operator $T$ is
{\bf near similar to uniform contraction} , iff there exists a bounded
operator $ Y\;>\;0 $     such that
$$
    TYT^*\;=\;\frac{Y}{I+Y}
$$

      We shall say, T is {\bf s-approximately similar to an unitary}  , iff 
there exists
 a net  $\{j_t\}_t$  of bounded uniformly positive operators such that

1) $1/M\,\le \,j_t\,\le M $ for some real $M$ (and every $t$ ),

2) $ strong-lim\,j_t\,=\,I $ (hence $strong-lim\,j_t^{-1}\,=\,I $),

3) for every fixed $t$ the operator $j_tT$ is similar to an unitary operator.
  $ \quad \Box$ 

\end{Definition}

 With this definition one can resume the section as follows:

\begin{Theorem}{ 1.} 
There exists an orthogonal decomposition 
$$
   H\;=\;H_<+\;H_=+\;H_> 
$$ 
such that

1) $H_<$ is $V$- surinvariant and the restriction of $V$
onto $H_<$ is near similar to an uniform contraction;

2) $H_>$ is $V^{*-1}$ - surinvariant and the restriction of $V^{*-1}$
onto $H_<$ is near similar to an uniform contraction;

3) If $H_<=\{0\}=H_>$
  then  $V$ is
s-aproximately similar to an unitary;

4)  $\;H_<\; \subset\; \overline{l_2(V)}\;$, 
$\quad H_>\;\subset\;\overline{l_2(V^{*-1})}$.

\end{Theorem}
 
\begin{Proof}. Set $H_<:=\overline{Ran\,Y_0}\, ,
                        H_=:=Ker\,X_0\cap Ker\,Y_0\, ,
                        H_>:=\overline{Ran\,X_0}      $

and apply the text.

\end{Proof}

\newpage

\section{ Stability of $U$ in terms of $Q_t$}

\begin{Proposition}{ 1.}
Let
$$
\frac{1}{M}I \le Q_t \le MI
$$
for some number $M > 0$ (and every $t$). 
Then $V$ is similar to an unitary operator. 
\end{Proposition}

\begin{Proof}. Let $Q_0$ be a weak limit point of the net
$Q_t, t \to +0$. It exists and
$$
\frac{1}{M}I \le Q_0 \le MI
$$
$$
Q_0 = V^*Q_0V
$$
(see Observation 2.2). Hence, $V$ is similar to an unitary operator.
\end{Proof}

\begin{Proposition}{  2.}
 Let $V$ be similar to an unitary operator.
Then
$$
\frac{1}{M}I \le Q_t \le MI
$$
for some number $M>0$ (and every $t$). 
\end{Proposition}

\begin{Proof}. For assumed $V$ there exists a number $M>0$
such that for every natural $n$ there hold
$$
\|V^n\|^2 \le M, \quad \|(V^*)^{-n}\| \le M
$$
Take arbitrary $x \in H$, number $z$, $0<z<1$ and apply the Observation 2.3:
$$
\begin{array}{rcl}
<x, Q_tx> & \le & (1-z) (\|Vx\|^2 +z\|V^2x\|^2 + 
\cdots +z^{n-1}\|V^nx\|^2+\cdots)
\\[0.5cm]
 & \le & (1-z){\displaystyle \frac{M}{1-z}} <x,x>=M<x,x>
\end{array}
$$
Hence, $Q_t \le MI$.

Now apply the Observation 2.3 to $Q_t^{-1}$ and $V^{*-1}$:
$$
\begin{array}{rcl}
<x,Q_t^{-1}x>&\le &(1-z)(\|V^{*-1}x \|^2+z\|(V^{*-1})^2x\|^2+ \cdots)\\[5pt]
              &\le & M<x,x>
\end{array}
$$
Hence, $Q_t^{-1}\le M$ and $\frac{1}{M}\le Q_t$.
\end{Proof}

\newpage
\section{ When  Spectrum has Dichotomy.}

Return us to the Example 1.1, which  was  called motivating .
We remarked there, that  $R_0$ %%%%is somewhat similar to
is reminiscent of one of the spectral projector of 
the  operator $V$   ( this  operator was taken there to be normal ).

Now we  will show that  somewhat similar  situation holds always,
especially  when the spectrum of the  operator $V$ 
does not  intersect the unit circle.

\begin{Observation}{ 1.} 
\fbox{$\displaystyle 
  Y_0=Y_0R_0= R_0Y_0 = R_0^{1/2}Y_0R_0^{1/2}
$}
In addition $R_0$  acts on $Ran\, Y_0$, hence on $\overline{ l_2(V) }$,
as identity operator
\end{Observation}

\begin{Proof}{\bf .} Recall 
$Y_t= tQ_t^{-1}$, 
$Y_0= s-\lim_{t\to +0} Y_t$,
$R_t=(I+Q_t)^{-1}$, 
$R_0 \in w-lim\ pt(R_t)$,
all these operators 
are selfadjoint. What is more, the straightforward calculation shows that
$$
Y_t (1-R_t) = \frac{t}{I+Q_t} = (1-R_t) Y_t  \,.
$$
The rest is obvious.
\end{Proof}

%It is straightforward to deduce now 
%\begin{Theorem 1.}

%\end{Theorem}

\begin{Theorem}{ 1.}
Suppose there is an 
 $V$-invariant subspace,
$L$ say,
such that
  $|spectrum(V|L)| < 1$; 
let $P$ denote  orthoprojector onto $L$ , 

Then
\begin{description}
\item[a)]
$$
 L \subset Ran\,Y_0^{1/2} \subset Ker(I-R_0)
$$

\item[b)]
$$
\|(I-R_t)P\|  \to 0  \quad (t\to +0)
$$
%\item[c)] 
%If in addition l_2 is closed then

\end{description}
\end{Theorem}
%%%%%%\end{document}
\begin{Proof}{\bf .}

{\bf a)}
Note $L\subset l_2(V)$ and apply theorem 2.1.

{\bf b)}

Return to the Observation 2.3:
$$
0\le Q_t\le (1-z)[V^*V+zV^{*2}V^2+\cdots+z^{n-1}V^{*n}V^n]+z^nV^{*n}Q_tV^n
$$
Let 
$V_P$  
denote 
$PV|L$.
Since
$VP=PVP$ , we can deduce that 
$$
0\leq PQ_tP
\leq 
(1-z)[V_P^*V_P+zV_P^{*2}V_P^2+
\cdots+z^{n-1}V_P^{*n}V_P^n]+z^nV_P^{*n}PQ_tPV_P^n
$$
%
%%%%%%%%% Next we will interpret the spectral condition argument 
%in accordance with standard practice:
Next we adopt the spectrum argument.
We make it 
in the same manner that standard practice suggests: 
there are some real positive  $\epsilon$, $M$  such that
$$
r(V_P)+\epsilon <1\,,\quad\|V_P^n\|  \leq M(r(V)+\epsilon)^n 
\qquad (n=1,2,\ldots) \,,
$$
Hence
$$\|z^nV_P^{*n}PQ_tPV_P^n\|\to 0 \quad (n\to \infty)   \,.$$
and it is routine matter to verify that
$$\|PQ_tP\| 
  \leq  \frac{(1-z)}{1-z(r(V_P)+\epsilon)^2}
                          M^2(r(V_P)+\epsilon)^2\,.$$
Note that  $M$ and $\epsilon$ does not depend on $z$ .
So, we obtain   $\|PQ_tP\|  \to 0 $ for $t\to +0$.

( One can show moreover : the  serie 
$$
(1-z)[V_P^*V_P+zV_P^{*2}V_P^2+\cdots+z^{n-1}V_P^{*n}V_P^n \cdots]
$$
is norm-convergent.)

   Now note  that 
$$
PR_tP\equiv 
P(I+Q_t)^{-1}P \geq P(I+PQ_tP)^{-1}P \stackrel{\|\ \|}{\to} 0
$$
$$
P(I-R_t)P\equiv 
PQ_t(I+Q_t)^{-1}P \leq PQ_tP(I+PQ_tP)^{-1} \stackrel{\|\ \|}{\to} 0
$$
Recall
$$0\leq I-R_t\leq I$$
Hence
$$
0\leq P(I-R_0)^2P
       =  P(I-R_0)^{1/2}(I-R_0)(I-R_0)^{1/2}P
            \leq P(I-R_0)P
$$
So, we can now establish that 
$$
\|(I-R_t)P\|^2  
  =\|P(I-R_t)(I-R_t)P\| 
  \leq \|P(I-R_t)P\| \to  0   \,,
$$
It was to be proved.
\end{Proof}

\begin{Corollary}{.}
Suppose that the spectrum of the  operator $V$ 
does not  intersect the unit circle;
let $P$ denote  orthoprojector onto spectral subspace $L$ 
 corresponded to the set 
$spectrum(V)\cap \{z\in {\bf C}| |z|<1\}$ .

Then
$$
   R_0=P \,.
$$

\end{Corollary}

\begin{Proof}{\bf .}
By Theorem 1
$$
(I-R_0)P=0  \,.
$$
Hence
$$
P = R_0 P   \,.
$$
Next note that the  equation
$$
  Q=V^*\frac{Q+t}{I+tQ} V,\quad Q\geq 0,\quad 0<t\leq 1 \qquad \eqno(^ *)
$$
is equivalent to the equation
$$
  Q^{-1}=V^{-1}\frac{Q^{-1}+t}{I+tQ^{-1}} V^{*-1},
\quad Q\geq 0,\quad 0<t\leq 1 \,.\qquad \eqno(^ **) 
$$
(for details see [Ch1,2])
\\
For a moment introduce for the (unique) solution of (*) 
a longer denotation: 
$$
 Q_t(V).
$$
It is straightforward to deduce now that
$$
Q_t(V^{*-1})=Q_t(V)^{-1}\,,\quad 
R_t(V)=I-R_t(V^{*-1})\,,\quad
Y_t(V)=X_t(V^{*-1})\, \ldots \mbox{etc.}
$$ 
Last recall the Standard Spectrum Theorems (see e.g. [RS])
and apply Theorem 1 to the operator 
$V^{*-1}$.
Then obtain
$$
R_0(I-P) =0  \,.
$$
Hence
$$
R_0=R_0 P  \,.
$$
To complete the proof let compare the second displayed formula 
with the last one in the current period.
\end{Proof}

\newpage
\bibliographystyle{unsrt}

\end{document}